# Magnetic Towers of Hanoi and their Optimal Solutions

*Uri Levy*

*Atlantium Technologies, Har-Tuv Industrial Park, Israel*

*uril@atlantium.com*

*August 5, 2010*

**Abstract**

The Magnetic Tower of Hanoi puzzle – a modified "base 3" version of the classical Tower of Hanoi puzzle as described in earlier papers, is actually a small set of independent sister-puzzles, depending on the "pre-coloring" combination of the tower's posts.

Starting with **R**ed facing up on a **S**ource post, working through an **I**ntermediate – colored or **N**eutral post, and ending **B**lue facing up on a **D**estination post, we identify the different pre-coloring combinations in (**S**,**I**,**D**) order. The Tower's pre-coloring combinations are {[(R,B,B) / (R,R,B)] ; [(R,B,N) / (N,R,B)] ; [(N,B,N) / (N,R,N)] ; [R,N,B] ; [(R,N,N) / (N,N,B)] ; [N,N,N]}.

In this paper we investigate these sister-puzzles, identify the algorithm that optimally solves each pre-colored puzzle, and prove its Optimality. As it turns out, five of the six algorithms, challenging on their own, are part of the algorithm solving the "natural", Free Magnetic Tower of Hanoi puzzle [N,N,N]. We start by showing that the N-disk Colored Tower [(R,B,B) / (R,R,B)] is solved by (3^N - 1)/2 moves.

Defining "Algorithm Duration" as the ratio of number of algorithm-moves solving the puzzle to the number of algorithm-moves solving the Colored Tower, we find the Duration-Limits for all sister-puzzles. In the order of the list above they are {[1] ; [10/11] ; [10/11] ; [8/11] ; [7/11] ; [20/33]}. Thus, the Duration-Limit of the Optimal Algorithm solving the Free Magnetic Tower of Hanoi puzzle is 20/33 or 606‰.

On the road to optimally solve this colorful Magnetic puzzle, we hit other "forward-moving" puzzle-solving algorithms. Overall we looked at 10 pairs of integer sequences. Of the twenty integer sequences, five are listed in the On-line Encyclopedia of Integer Sequences, the other fifteen – not yet.

The large set of different solutions is a clear indication to the freedom-of-wondering that makes this Magnetic Tower of Hanoi puzzle so colorful.

## 1. Introduction

The Magnetic Tower of Hanoi (MToH) puzzle is a modified version of the Classical Tower of Hanoi (ToH) puzzle. While the Classical ToH version spans base 2, the far more challenging MToH puzzle spans base 3. I described the MToH puzzle and analyzed its solutions in terms of puzzle-solving Algorithms and in terms of number of moves in two earlier papers. A first short version[1] and a second, revised and more complete version (reference [2] and references therein). However, in these two earlier versions I did NOT present the Optimal Algorithms (the Algorithms that solve the MToH puzzles with minimum number of moves).

Presenting these Optimal Algorithms and proving their Optimality is the objective of this third paper.

For the sake of brevity, and in view of the previously published papers, The MToH puzzle and its solving-rules are not repeated here. We can thus move on to the next section which is an overview of the MToH-solving Algorithms.

## 2. Overview of the MToH-solving Algorithms

The MToH puzzle is actually a set of sister-puzzles, depending on the combination of "pre-coloring" of the posts. We have identified six such combinations – Table 1.

Of the six pre-coloring combinations one {[NBN / NRN]} is equivalent to the {[RBN / NRB]} combination in terms of the solving Algorithm, and in addition, this particular combination – it turns out – does not participate in the solution of the Free-MToH. Thus, while the [NBN / NRN] combination is listed in Table 1 (and is depicted in Figure 1), it is not counted, and this particular combination is not discussed further in this paper.

The dark-green rows in Table 1 designate the Optimal Algorithms.

"Duration" in Table 1 is the limit for large number of disks of the ratio of number of algorithm-moves solving the puzzle to the number of algorithm-moves solving the Colored Tower.

The number in the Algorithm designation, in the Disk-move series designation and in the total-move series designation represents approximation to the solution's Duration-Limit in percent (two digits) or in promil (three and four digits).

"OEIS" In Table 1 stands for the On-line Encyclopedia of Integer Sequences[3].

Note that RRN and NBB combinations imply RRB and RBB respectively and are thus represented by the Colored MToH combinations listed in row number 1.

A name ("Colored", "Free", etc.) is given to each sister of the MToH puzzle family.

| Alg. | MToH | Pre-coloring | Algorithm | Disk-move | OEIS | Total-move | OEIS | Duration |
|---|---|---|---|---|---|---|---|---|
| # | State | S ; I ; D | Designation | Series | Y / N | Series | Y / N | Limit |
| 1 | Colored | RRB / RBB | RRB1000 | P1000(k) | YES | S1000(N) | YES | 1 |
| 2 | Semi-Free-C | RBN / NRB | RBN909 | P909(k) | NO | S909(N) | NO | 10/11 |
| X | Nearly-Free-N | NBN / NRN | NBN909 | P909(k) | X | S909(N) | X | 10/11 |
| 3.1 | Semi-Free-C | RNB | RNB75 | P75(k) | YES | S75(N) | NO | 3/4 |
| 3.2 | Semi-Free-C | RNB | RNB727 | P727(k) | NO | S727(N) | NO | 8/11 |
| 3.3 | Nearly-Free-C | RNN / NNB | RNN67 | P67(k) | YES | S67(N) | YES | 2/3 |
| 4.1 | Nearly-Free-C | RNN / NNB | RNN64 | P64(k) | NO | S64(N) | NO | 23/36 |
| 4.2 | Nearly-Free-C | RNN / NNB | RNN636 | P636(k) | NO | S636(N) | NO | 7/11 |
| 5.1 | Free | NNN | NNN62 | P62(k) | NO | S62(N) | NO | 67/108 |
| 5.2 | Free | NNN | NNN61 | P61(k) | NO | S61(N) | NO | 197/324 |
| 5.3 | Free | NNN | NNN606 | P606(k) | NO | S606(N) | NO | 20/33 |

Table 1: *The six pre-colored combinations of the MToH puzzle. The [NBN / NRN] combination is not numbered because it is equivalent to the [RBN / NRB] combination and because it does not participate in the Free-MToH puzzle solution. Paired combinations are "Time-Reversal Pairs" and are necessarily solved by "similar" Algorithms (obey the same recurrence relations and generate a single integer sequence).*
*The dark-green rows designate the Optimal Algorithms. The number in the Algorithm designation, in the Disk-move series designation and in the Total-move series designation represents approximation to the solution's Duration-Limit in percent (two digits) or in promil (three and four digits).*
*OEIS stands for the On-line Encyclopedia of Integer Sequences.*
*Note that RRN and NBB imply RRB and RBB and are thus represented by the Colored MToH combinations (row number 1).*

An explicit pictorial description of the six pre-colored combinations of the MToH sister-puzzles is shown in Figure 1. The one four-digit number and the rest three-digit numbers in the table, represent an approximation to the Duration-Limit of the MToH **Optimal** Solution (in promil).

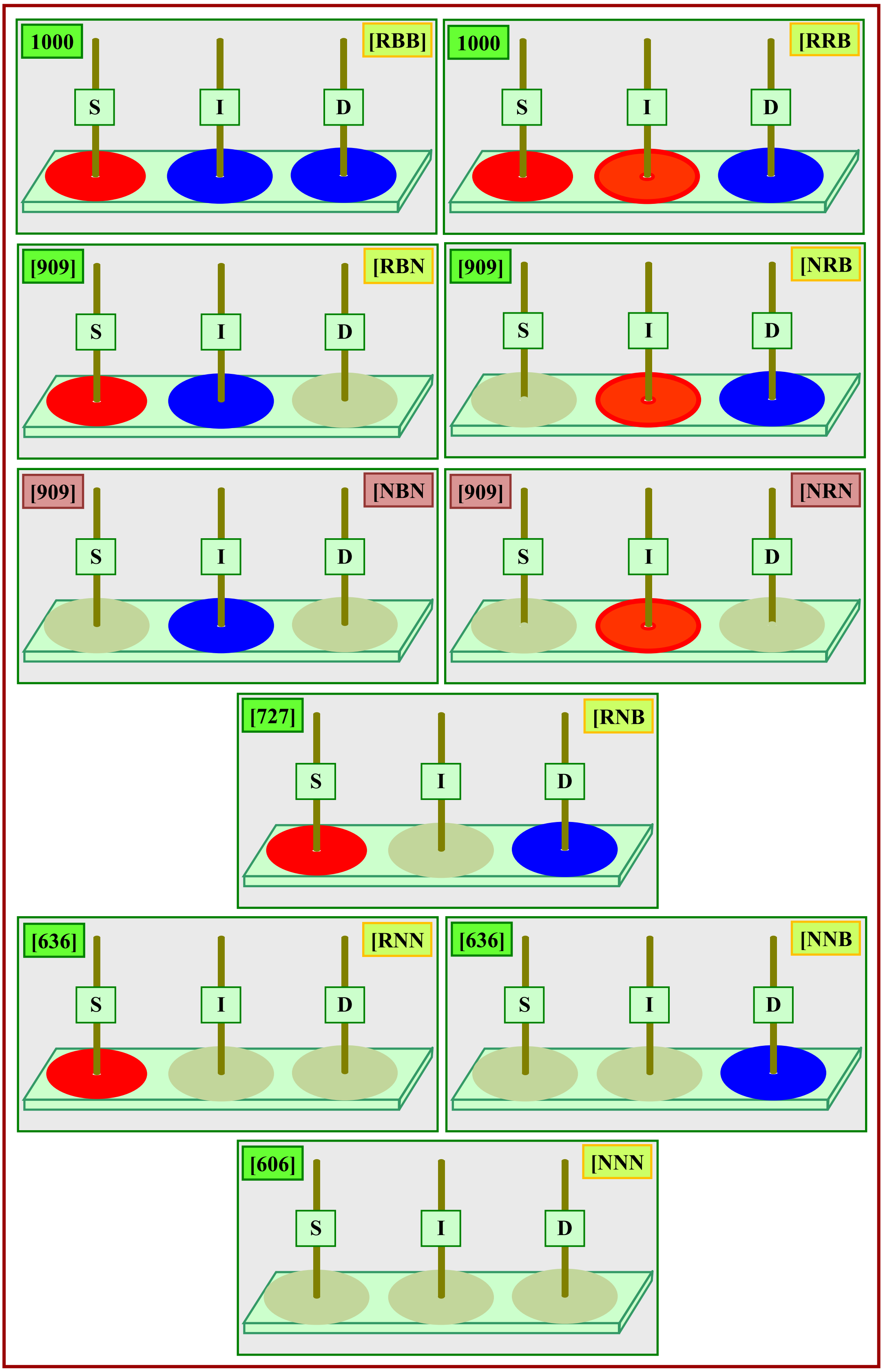

Figure 1: *Pre-colored combinations of the MToH sister-puzzles. The numbers in the top-left green box represent the Duration-Limit of the Optimal solution.*

Before moving on to discussing the Optimal Solution-Algorithms, let's just list the integer sequences generated by the "forward-moving" yet non-Optimal Solution-Algorithms.

### 2.1. Integer Sequences generated by the non-Optimal Solution-Algorithms

For paper completeness, integer sequences generated by non-Optimal forward-moving MToH puzzle-Solutions are listed in Table 2 and Table 3. The first sequence in each table (RRB1000 – dark green) is of course Optimal and is listed as a reference. The RRB1000 Algorithm is analyzed in detail in section 3 below.

Table 2 lists integer sequences of disk-moves - the number of moves each disk makes during execution of the particular Algorithm, given the total number of disks in the stack. Disk numbering is from bottom to top – largest disk's number is k = 1 and smallest disk's number is k = N (N = 20 for Table 2). Designation of each disk-move sequence generated by an xy Algorithm is Pxy(k).

Table 3 lists the sequences of total number of moves executed to solve the MToH puzzle by the particular Algorithm. Designation of the total number of moves sequence generated by an xy Algorithm is Sxy(N).

Also presented in Table 2 and Table 3 are **closed-form expressions** developed for each tabulated sequence.

Duration Limits, both numerically calculated (for k = 20 or N = 20) and exact (deduced from the closed form expressions) are listed for each Algorithm as well.

Finally – OEIS appearance of the integer sequence is noted at the last row of each table.

| | Colored | SemiFree-C | NearlyFree-C | NearlyFree-C | Free | Free |
|---|---|---|---|---|---|---|
| | "1000" | "75" | "67" | "64" | "62" | "61" |
| | RRB1000 | RNB75 | RNN67 | RNN64 | NNN62 | NNN61 |
| Y | 1 | 3/4 | 2/3 | 23/36 | 67/108 | 197/324 |
| k-Odd<br>k-Even | 3^(k-1) {≡Z} | Y*Z+1/4<br>Y*Z+3/4 | Y*Z+1 | Y*Z+k-7/4<br>Y*Z+k-9/4 | Y*Z+11/4<br>Y*Z+9/4 | Y*Z+2*k-25/4<br>Y*Z+2*k-27/4 |
| | k > 0 | k > 0 | k > 1 | k > 2 | k > 3 | k > 4 |
| K | $P_{1000}(k)$ | $P_{75}(k)$ | $P_{67}(k)$ | $P_{64}(k)$ | $P_{62}(k)$ | $P_{61}(k)$ |
| 1 | 1 | 1 | 1 | 1 | 1 | 1 |
| 2 | 3 | 3 | 3 | 3 | 3 | 3 |
| 3 | 9 | 7 | 7 | 7 | 7 | 7 |
| 4 | 27 | 21 | 19 | 19 | 19 | 19 |
| 5 | 81 | 61 | 55 | 55 | 53 | 53 |
| 6 | 243 | 183 | 163 | 159 | 153 | 153 |
| 7 | 729 | 547 | 487 | 471 | 455 | 451 |
| 8 | 2187 | 1641 | 1459 | 1403 | 1359 | 1339 |
| 9 | 6561 | 4921 | 4375 | 4199 | 4073 | 4001 |
| 10 | 19683 | 14763 | 13123 | 12583 | 12213 | 11981 |
| 11 | 59049 | 44287 | 39367 | 37735 | 36635 | 35919 |
| 12 | 177147 | 132861 | 118099 | 113187 | 109899 | 107727 |
| 13 | 531441 | 398581 | 354295 | 339543 | 329693 | 323149 |
| 14 | 1594323 | 1195743 | 1062883 | 1018607 | 989073 | 969409 |
| 15 | 4782969 | 3587227 | 3188647 | 3055799 | 2967215 | 2908187 |
| 16 | 14348907 | 10761681 | 9565939 | 9167371 | 8901639 | 8724515 |
| 17 | 43046721 | 32285041 | 28697815 | 27502087 | 26704913 | 26173497 |
| 18 | 129140163 | 96855123 | 86093443 | 82506231 | 80114733 | 78520437 |
| 19 | 387420489 | 290565367 | 258280327 | 247518663 | 240344195 | 235561255 |
| 20 | 1162261467 | 871696101 | 774840979 | 742555955 | 721032579 | 706683703 |
| T(20) | 1 | 0.750000001 | 0.666666668 | 0.638888904 | 0.620370372 | 0.60802472 |
| T_limit | 1 | 3/4 | 2/3 | 23/36 | 67/108 | 197/324 |
| OEIS | YES | YES | YES | NO | NO | NO |

Table 2: *First twenty elements of the disk-moves of integer sequences generated by the non-Optimal Algorithms (except for the RRB1000 sequence in the dark green column on the left).*

| | Colored | SemiFree-C | NearlyFreeC | NearlyFree-C | Free | Free |
|---|---|---|---|---|---|---|
| | "1000" | "75" | "67" | "64" | "62" | "61" |
| | RRB1000 | RNB75 | RNN67 | RNN64 | NNN62 | NNN61 |
| Y | 1 | 3/4 | 2/3 | 23/36 | 67/108 | 197/324 |
| N-Odd | Z - 1/2 {Z ≡ (3^N)/2} | Y*Z+0.5*N-5/8 | Y*Z+N-1 | Y*Z+0.5*N^2-1.5*N+19/8 | Y*Z+2.5*N-39/8 | Y*Z+N^2-5.5*N+93/8 |
| N-Even | | Y*Z+0.5*N-3/8 | | Y*Z+0.5*N^2-1.5*N+17/8 | Y*Z+2.5*N-41/8 | Y*Z+N^2-5.5*N+91/8 |
| | N > 0 | N > 0 | N > 0 | N > 1 | N > 2 | N > 3 |
| N | $S_{1000}(N)$ | $S_{75}(N)$ | $S_{67}(N)$ | $S_{64}(N)$ | $S_{62}(N)$ | $S_{61}(N)$ |
| 1 | 1 | 1 | 1 | 1 | 1 | 1 |
| 2 | 4 | 4 | 4 | 4 | 4 | 4 |
| 3 | 13 | 11 | 11 | 11 | 11 | 11 |
| 4 | 40 | 32 | 30 | 30 | 30 | 30 |
| 5 | 121 | 93 | 85 | 85 | 83 | 83 |
| 6 | 364 | 276 | 248 | 244 | 236 | 236 |
| 7 | 1093 | 823 | 735 | 715 | 691 | 687 |
| 8 | 3280 | 2464 | 2194 | 2118 | 2050 | 2026 |
| 9 | 9841 | 7385 | 6569 | 6317 | 6123 | 6027 |
| 10 | 29524 | 22148 | 19692 | 18900 | 18336 | 18008 |
| 11 | 88573 | 66435 | 59059 | 56635 | 54971 | 53927 |
| 12 | 265720 | 199296 | 177158 | 169822 | 164870 | 161654 |
| 13 | 797161 | 597877 | 531453 | 509365 | 494563 | 484803 |
| 14 | 2391484 | 1793620 | 1594336 | 1527972 | 1483636 | 1454212 |
| 15 | 7174453 | 5380847 | 4782983 | 4583771 | 4450851 | 4362399 |
| 16 | 21523360 | 16142528 | 14348922 | 13751142 | 13352490 | 13086914 |
| 17 | 64570081 | 48427569 | 43046737 | 41253229 | 40057403 | 39260411 |
| 18 | 193710244 | 145282692 | 129140180 | 123759460 | 120172136 | 117780848 |
| 19 | 581130733 | 435848059 | 387420507 | 371278123 | 360516331 | 353342103 |
| 20 | 1743392200 | 1307544160 | 1162261486 | 1113834078 | 1081548910 | 1060025806 |
| T(20) | 1 | 0.750000006 | 0.666666678 | 0.638888988 | 0.620370396 | 0.608024864 |
| T_limit | 1 | 3/4 | 2/3 | 23/36 | 67/108 | 197/324 |
| OEIS | YES | NO | YES | NO | NO | NO |

Table 3: *First twenty elements of the total-moves integer sequence generated by the non-Optimal Algorithms (except for the RRB1000 sequence in the dark green column on the left).*

### 2.2. Integer Sequences generated by the Optimal Solution-Algorithms

As part of the Overview of the MToH-solving Algorithms, integer sequences generated by the Optimal Solution-Algorithms are listed in Table 4 and in Table 5. Closed-form expressions are not included in these tables (they are presented in subsequent sections).

The rest of the paper is devoted to a detailed discussion of each of these Optimal Solution-Algorithms.

| | Colored | SemiFree-N | SemiFree-C | NearlyFree | Free |
|---|---|---|---|---|---|
| | "1000" | "909" | "727" | "636" | "606" |
| | RRB1000 | RBN909 | RNB727 | RNN636 | NNN606 |
| K | $P_{1000}(k)$ | $P_{909}(k)$ | $P_{727}(k)$ | $P_{636}(k)$ | $P_{606}(k)$ |
| 1 | 1 | 1 | 1 | 1 | 1 |
| 2 | 3 | 3 | 3 | 3 | 3 |
| 3 | 9 | 9 | 7 | 7 | 7 |
| 4 | 27 | 25 | 21 | 19 | 19 |
| 5 | 81 | 75 | 61 | 55 | 53 |
| 6 | 243 | 223 | 179 | 159 | 153 |
| 7 | 729 | 665 | 535 | 471 | 451 |
| 8 | 2187 | 1993 | 1597 | 1403 | 1339 |
| 9 | 6561 | 5971 | 4781 | 4191 | 3997 |
| 10 | 19683 | 17903 | 14331 | 12551 | 11961 |
| 11 | 59049 | 53697 | 42967 | 37615 | 35835 |
| 12 | 177147 | 161065 | 128869 | 112787 | 107435 |
| 13 | 531441 | 483163 | 386557 | 338279 | 322197 |
| 14 | 1594323 | 1449439 | 1159587 | 1014703 | 966425 |
| 15 | 4782969 | 4348233 | 3478647 | 3043911 | 2899027 |
| 16 | 14348907 | 13044585 | 10435757 | 9131435 | 8696699 |
| 17 | 43046721 | 39133571 | 31306989 | 27393839 | 26089517 |
| 18 | 129140163 | 117400431 | 93920555 | 82180823 | 78267673 |
| 19 | 387420489 | 352200881 | 281761015 | 246541407 | 234801675 |
| 20 | 1162261467 | 1056601993 | 845282069 | 739622595 | 704402987 |
| T(20) | 1 | 0.909091476 | 0.727273589 | 0.636365066 | 0.606062411 |
| **T-limit** | 1 | 10/11 | 8/11 | 7/11 | 20/33 |
| OEIS | YES | NO | NO | NO | NO |

Table 4: *First twenty elements of the disk-moves of integer sequences generated by the Optimal Algorithms.*

| | Colored | SemiFree-N | SemiFree-C | NearlyFree | Free |
|---|---|---|---|---|---|
| | "1000" | "909" | "727" | "636" | "606" |
| | RRB1000 | RBN909 | RNB727 | RNN636 | NNN606 |
| N | $S_{1000}(N)$ | $S_{909}(N)$ | $S_{727}(N)$ | $S_{636}(N)$ | $S_{606}(N)$ |
| 1 | 1 | 1 | 1 | 1 | 1 |
| 2 | 4 | 4 | 4 | 4 | 4 |
| 3 | 13 | 13 | 11 | 11 | 11 |
| 4 | 40 | 38 | 32 | 30 | 30 |
| 5 | 121 | 113 | 93 | 85 | 83 |
| 6 | 364 | 336 | 272 | 244 | 236 |
| 7 | 1093 | 1001 | 807 | 715 | 687 |
| 8 | 3280 | 2994 | 2404 | 2118 | 2026 |
| 9 | 9841 | 8965 | 7185 | 6309 | 6023 |
| 10 | 29524 | 26868 | 21516 | 18860 | 17984 |
| 11 | 88573 | 80565 | 64483 | 56475 | 53819 |
| 12 | 265720 | 241630 | 193352 | 169262 | 161254 |
| 13 | 797161 | 724793 | 579909 | 507541 | 483451 |
| 14 | 2391484 | 2174232 | 1739496 | 1522244 | 1449876 |
| 15 | 7174453 | 6522465 | 5218143 | 4566155 | 4348903 |
| 16 | 21523360 | 19567050 | 15653900 | 13697590 | 13045602 |
| 17 | 64570081 | 58700621 | 46960889 | 41091429 | 39135119 |
| 18 | 193710244 | 176101052 | 140881444 | 123272252 | 117402792 |
| 19 | 581130733 | 528301933 | 422642459 | 369813659 | 352204467 |
| 20 | 1743392200 | 1584903926 | 1267924528 | 1109436254 | 1056607454 |
| T(20) | 1 | 0.909092014 | 0.727274407 | 0.636366421 | 0.606064117 |
| **T-limit** | 1 | 10/11 | 8/11 | 7/11 | 20/33 |
| OEIS | YES | NO | NO | NO | NO |

Table 5: *First twenty elements of the total-moves integer sequences generated by the Optimal Algorithms.*

Last in this Algorithm overview section are Duration curves for the Optimal Solution-Algorithms.

### 2.3. Duration curves for the Optimal MToH Solution-Algorithms

Figure 2, as part of the MToH Solution-Algorithm overview, presents Duration curves for the Optimal Algorithms. "Duration" for Algorithm xyz solving an N-disk puzzle [ $T_{xyz}(N)$ ] is defined as

$$T_{xyz}(N) \equiv S_{xyz}(N) / S_{1000}(N) . \qquad (1)$$

Clearly, a less color-restrictive MToH puzzle is solved (efficiently) in a smaller number of moves.

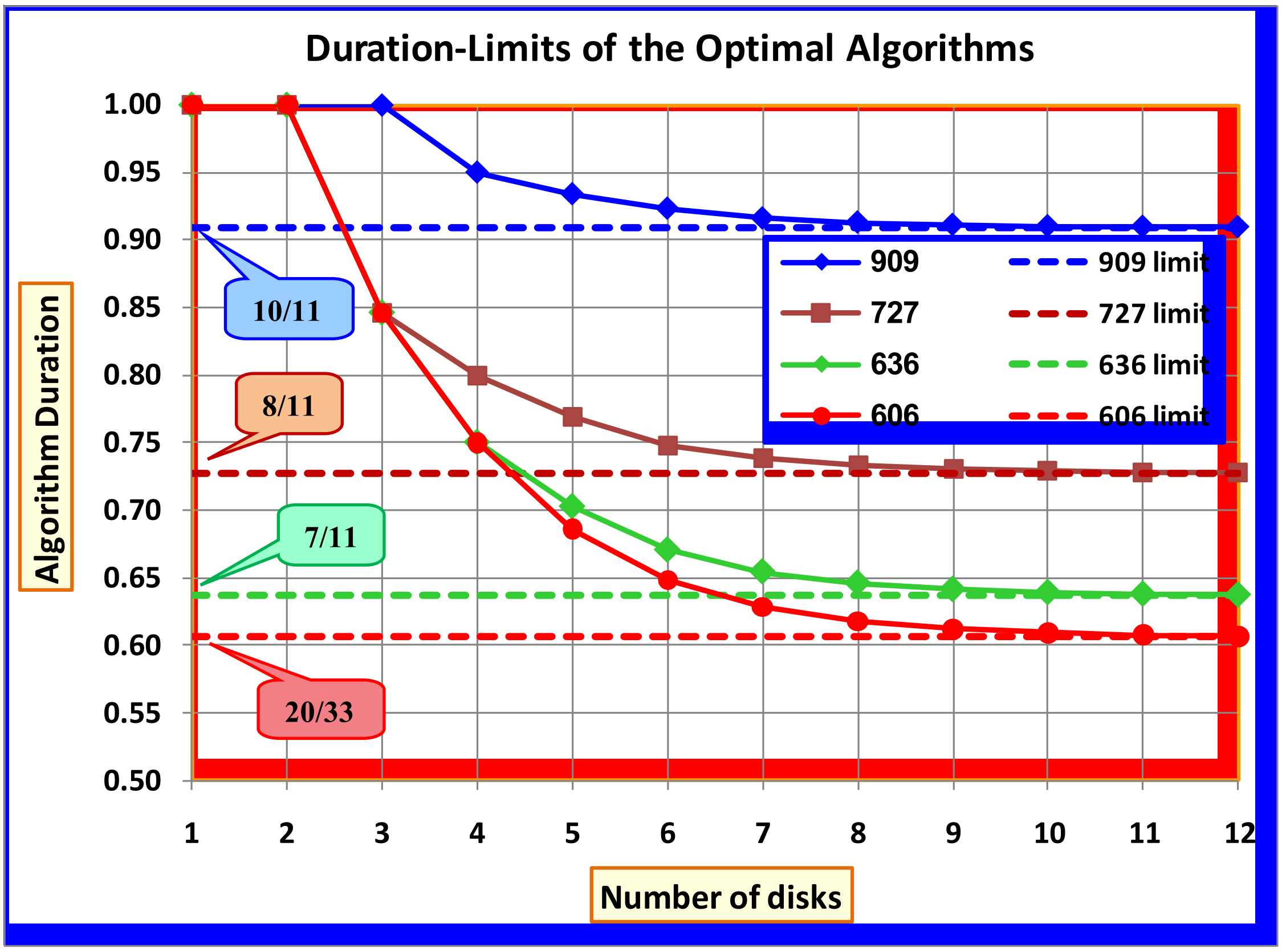

Figure 2: *Duration curves for the Optimal Algorithms. These Algorithms are discussed in detail in section 3 below.*

So much for the overview of the MToH puzzle-solving Algorithms. Let's take a closer look now at the Optimal Solution-Algorithms.

## 3. The Optimal Solution-Algorithms

In this section, the five Optimal Solution-Algorithms are discussed in detail. In Table 1, these Optimal Solution-Algorithms are numbered {1 ; 2 ; 3.2 ; 4.2 ; 5.3}.

For all five Algorithms, the solving task calls for moving N **RED** facing up disks orderly stacked on a **S**ource post (**S**) to a **BLUE** facing up disks orderly stacked on a **D**estination post (**D**), using an **I**ntermediate post (**I**). A move consists of transporting a single **flipped** disk from one post to another, obeying the two MToH move rules – the "Size-Rule" and the "Magnetic-Rule"[2].

As an example, an MToH start-state (number 3.2 in Table 1 – the Semi-Free-C combination, arbitrarily selected), is shown in Figure 3.

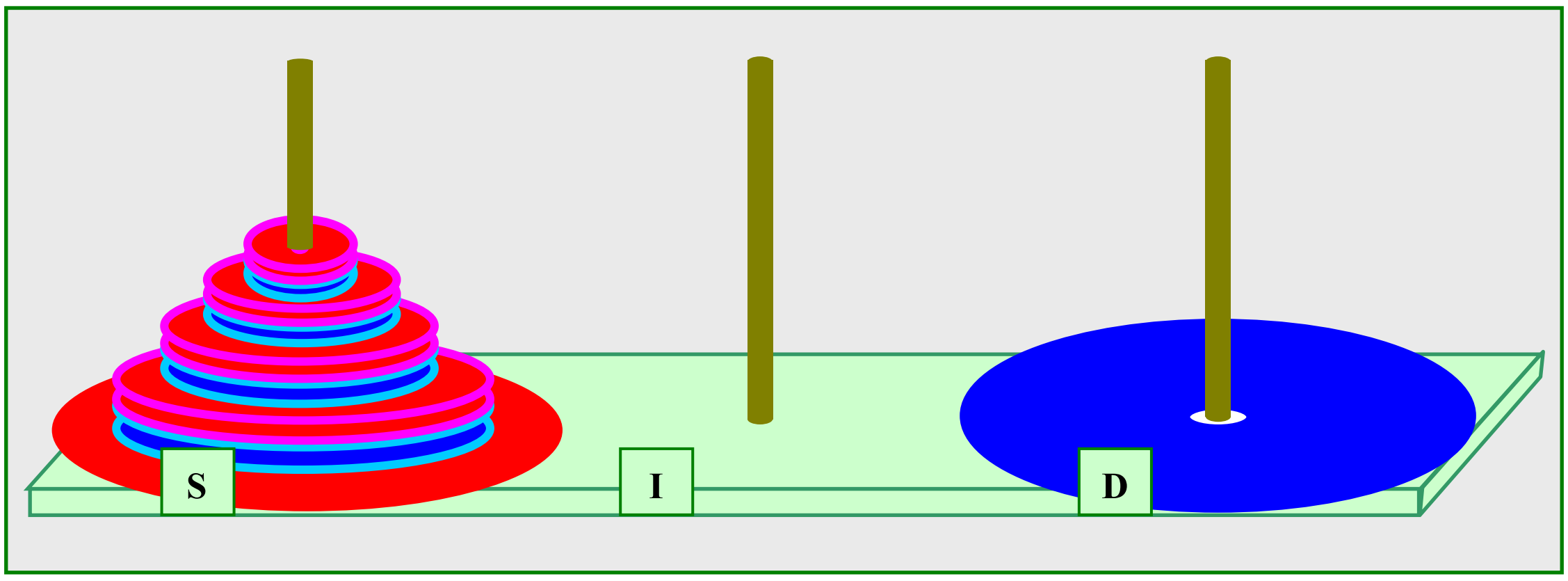

Figure 3: *An example of the MToH start-state with an arbitrarily selected pre-coloring configuration (the Semi-Free-C MToH puzzle in this case). The puzzle-solving task calls for moving N **RED** facing up disks orderly stacked on a **S**ource post (**S**) to **BLUE** facing up disks orderly stacked on a **D**estination post (**D**), using an **I**ntermediate post (**I**).*

Note: The "**BLUE** facing up" requirement for the MToH end-state can be replaced by a "minimum number of moves" requirement, since a **RED**-to-**RED** MToH puzzle is always (for any number of disks and for any pre-coloring configuration) solved with greater number of moves. This statement is not re-visited below and is not proved in this paper.

We start the Optimal-Algorithms discussion with number 1 in Table 1 – the Algorithms solving the Colored MToH.

### 3.1. The Colored MToH and its solving Optimal Algorithms

Let's start by reminding ourselves of the pre-coloring configuration for the Colored MToH-puzzle.

#### 3.1.1. The pre-coloring configuration of the Colored MToH

The pre-coloring configuration for the Colored MToH puzzle is shown in Figure 4.

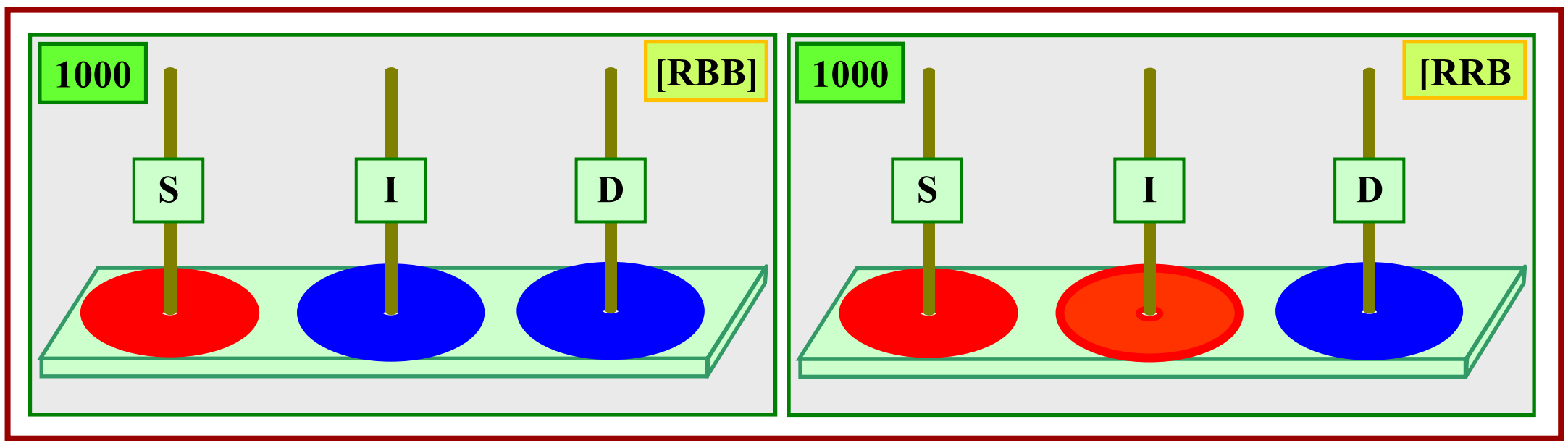

Figure 4: *The pre-coloring configurations of the Colored MToH puzzle. The two distinct puzzle-solving Algorithms form a Time-Reversal Pair.*

Next – the Optimal "1000" puzzle-solving Algorithms.

#### 3.1.2. The RBB1000 / RRB1000 Optimal Algorithms

The Optimal Algorithms solving the Colored MToH puzzle are listed in Table 6.

| | SRBB1000 | SRRB1000 |
|---|---|---|
| | function move_SRBB1000(n,s,d,i) | function move_SRRB1000(n,s,d,i) |
| | j = %N_disks + 1 - n | j = %N_disks + 1 - n |
| | if n > 0 | if n > 0 |
| 1 | move_SRBB1000(n-1,s,i,d) | move_SRRB1000(n-1,s,d,i) |
| 2 | move(j,s,d) | move_SRBB1000(n-1,d,i,s) |
| 3 | move_SRRB1000(n-1,i,s,d) | move(j,s,d) |
| 4 | move_SRBB1000(n-1,s,d,i) | move_SRRB1000(n-1,i,d,s) |

Table 6: *The Optimal Algorithms solving the Colored MToH-puzzle.*

The RBB pre-colored configuration is solved by the SRBB1000 Algorithm listed on the left of Table 6, while the RRB pre-colored configuration is solved by the SRRB1000 Algorithm listed on the right of Table 6.

Every MToH solving Algorithm generates a pair of recurrence relations – one for the total number of moves (designated for example S1000(N) – see Table 5) and one for the number of disk-moves (designated for example P1000(k) – see Table 4). The "S" in each of the function calls in Table 6 is introduced with this "S" vs. "P" recursive relations difference in mind (see subsection 3.1.5. below).

The " j = %N_disks + 1 – n" command introduced in each of the listed functions in Table 6 is a transformation to accommodate the disk numbering convention. Thus - move(j,s,d) is actually move(1,s,d) when n = N_disks [and move(j+1,s,d), as appears in other Algorithms, is actually move(2,s,d) when n = N_disks].

The SRBB1000 Algorithm and the SRRB1000 Algorithm form a Time-Reversal Pair as explained in the next section.

### 3.1.3. A Time-Reversal Algorithm Pair

When solving the MToH puzzle "forward", it is **always possible** to go backwards, in a **T**ime-**R**eversal fashion, and solve the MToH puzzle this way, from the end to the beginning. If the **RED**-**BLUE** colors are now swapped and Destination and Source posts are swapped, the Time-Reversed backward Algorithm becomes a "legitimate" forward solving Algorithm.

This Time-Reversal operation is easily executed and fully appreciated when playing an MToH-puzzle-applet, which allows multiple "undo" operations (now on-line[5],[6]).

Clearly, the forward Algorithm and the backward Algorithm form a Time-Reversal Algorithm Pair. We refer to the two Algorithms forming the pair as "Brothers".

The two Algorithms listed in Table 6 form such a Time-Reversal Algorithm Pair:

$$TR(\text{SRBB1000}) \equiv \text{SRRB1000} \quad and \quad TR(\text{SRRB1000}) \equiv \text{SRBB1000} \quad ,(2)$$

where "*TR*" signifies a Time-Reversal operation.

Proof of the statement in Equation 2 is left for the reader.

Not less clear is the fact that the total number of puzzle-solving moves found for one member of a Time-Reversal Algorithm Pair, equals the total number of puzzle-solving moves found for its Brother, for any stack-height - N. And similarly for the number of moves of the individual disks. These equalities hold of course for any Time-Reversal Algorithm Pair – Optimal or not.

Table 7 lists explicitly the sequence of number of disks on each post as the three-disk Colored MToH puzzle is solved by the two Optimal Algorithms in question. Inspecting the two columns reveals the Time-Reversal nature of these two Algorithm-Brothers.

| Move # | SRBB1000 | SRRB1000 |
|---|---|---|
| 0 | 3 0 0 | 3 0 0 |
| 1 | 2 0 1 | 2 0 1 |
| 2 | 1 1 1 | 2 1 0 |
| 3 | 2 1 0 | 1 1 1 |
| 4 | 1 2 0 | 1 0 2 |
| 5 | 0 2 1 | 2 0 1 |
| 6 | 1 1 1 | 2 1 0 |
| 7 | 0 1 2 | 1 1 1 |
| 8 | 1 0 2 | 1 2 0 |
| 9 | 2 0 1 | 0 2 1 |
| 10 | 1 1 1 | 0 1 2 |
| 11 | 0 1 2 | 1 1 1 |
| 12 | 1 0 2 | 1 0 2 |
| 13 | 0 0 3 | 0 0 3 |

Table 7: *Number of disks on each post when solving the three-disk Colored MToH-puzzle by each of the Optimal Algorithms. The table does NOT specify disk size, but the reader can appreciate that the (111) state on the left in rows 2, 6, 10 each represents a* ***unique*** *MToH state due to its unique disk-size arrangement. Posts' colors are not stated either, but in this specific case post colors are known and are* ***fixed*** *throughout the solving procedure.*

*Reading the three digit numbers in one column from right to left and from bottom to top, they are seen to be equal to the ordinary numbers on the other column read from top to bottom.*

Back to Optimality, we will now prove that each of the two Colored MToH Solving-Algorithms listed in Table 6 is Optimal.

### 3.1.4. Proof of Optimality for the "1000" Pair

The Optimality proof for both Algorithms in Table 6 is a coupled-recursive proof. We show that both Algorithms are Optimal for N = 1, we assume that **both** Algorithms are Optimal for any n = N, and we prove sequentially that each Algorithm is Optimal for n = N + 1. The N + 1 part of the proof is based on step-by-step inspection of the MToH puzzle and using a "must" argument for each step.

Let's see.

We start with the SRBB1000 Algorithm, on the left of Table 6.

For N = 1, both Algorithms call for one move to solve the Colored MToH puzzle. "One move" Algorithm is obviously Optimal.

For n = N – 1 we **assume** Optimality of both Algorithms (N > 1).

For n = N:

The first step (line 1 in Table 6) is applying the SRBB1000 Algorithm itself to transport N-1 disks (disk 2 to disk N) from **RED-S** to **BLUE-I** using **D**. This step is Optimal by assumption. The presence of the big disk at the bottom of post **S** (was not present there in the N-1 case) makes no difference because it is big (biggest) and because it is **RED**. And we **must** perform this (Optimal) step in order to free the big disk laying on the **S**ource post and we must clear the **D**estination post.

Note: The designations of the upper case (**S**,**I**,**D**) letters used in the text and the designations of the lower case (**s**,**i**,**d**) letters used in the functions of Table 6 are identical. Both are short for **S**ource, **I**ntermediate, **D**estination (posts).

The next step (line 2 in Table 6) is moving the big disk (disk number 1) from **S** to **D**. We **must** move any disk at least once. In this case (and in all other cases) we move the big disk exactly once which is certainly Optimal.

The third step (line 3 in Table 6) is moving N-1 disks from **I** to **S**. We **must** move all N-1 disks to the **RED**-colored Source post because we must move disk number 2 to the **RED**-colored Source post and the smaller disks, after parking on the **BLUE**-colored Destination post must all fold back on the **RED**-colored Source post. But now, for the task in question, the "pre-coloring" state of the Tower [in the (**S**,**I**,**D**) order] is BBR which is equivalent to RRB, so we **have** to use the Time-Reversal Brother Algorithm (SRRB1000 – on the right of Table 6) to execute this step. Here again, the presence of the big disk at the bottom of the Destination post (was not present in the N-1 case) makes no difference because it is big (biggest) and because it is (necessarily) BLUE. The N-1 disk transport by the SRRB1000 Algorithm is Optimal by assumption.

Finally (line 4 in Table 6), we **must** move all N-1 disks from **S** to **D**. The pre-coloring state (for the task) is again RBB so we resort again to the services of the original SRBB1000 Algorithm. The transport, for this step too, is not affected by the presence of the big disk on the Destination post and it is Optimal by assumption.

Puzzle solved. And the (conditional) proof of Optimality of the SRBB1000 Algorithm ends here.

For the Time-Reversal Brother Algorithm - SRRB1000, we follow the same path of proof, using similar arguments. Note that even now, when we use SRBB1000 (line 2 on the right of Table 6), we still **assume** that it is Optimal because the proof is coupled and it is incomplete until this second part of SRRB1000 Optimality proof ends.

But now, when we are done with the second Brother, we **know** that both Algorithms are Optimal.

End of Optimality proof for the two "1000" Brothers.

In subsequent proofs, when we run into one of these "1000" pre-coloring configurations (happens rather frequently) and we execute a step by one or the other of these "1000" Time-Reversal Pair Algorithms, we know it is Optimal.

Next – on to recurrence relations.

### 3.1.5. Recurrence relations for the "1000" Pair

Given the "1000" Solution-Algorithm Brothers (Table 6), we can extract recurrence relations for the associated number of moves.

First, necessarily for any Time-Reversal Algorithm Pair –

$$S_{SRBB1000}(N) = S_{SRRB1000}(N) \equiv S_{1000}(N) \tag{3A}$$

$$P_{SRBB1000}(k) = P_{SRRB1000}(k) \equiv P_{1000}(k) \tag{3B}$$

where $S_{xyz}(N)$ is defined next to Equation 1 and $P_{xyz}(k)$ is the number of moves of disk number k during solving the MToH puzzle by Algorithm xyz (independent of the total number (N) of disks in the stack).

From the left part of Table 6 and using Equation 3A:

$$S_{SRBB1000}(N+1) = 2 \cdot S_{SRBB1000}(N) + S_{SRRB1000}(N) + 1 = 3 \cdot S_{SRBB1000}(N) + 1 \tag{4A}$$

In general, the recurrence relations for the total number of moves (Equation 4A in this particular case) **must** work for the recurrence relations of the moves of any disk (k), only without the "singles" ("1" in this particular case – Equation 4A). The singles are not counted because they always apply only to the big disk (or disks) at the bottom of the stack. So we have –

$$P_{SRBB1000}(k+1) = 2 \cdot P_{SRBB1000}(k) + P_{SRRB1000}(k) = 3 \cdot P_{SRBB1000}(k) \tag{4B}$$

And using Equation 3 again, we can finally write:

$$S_{1000}(N+1) = 3 \cdot S_{1000}(N) + 1 \quad ; \quad S_{1000}(1) = 1 \tag{5A}$$

$$P_{1000}(k+1) = 3 \cdot P_{1000}(k) \quad ; \quad P_{1000}(1) = 1 \tag{5B}$$

The Recurrence Relations 5A and 5B hold of course for both "1000" Algorithms.

#### 3.1.6. Closed-form expressions for the "1000" Algorithm

It is not too difficult to show (prove) that the Recurrence Relations 5A and 5B hold **if and only if** they generate the following closed-form expressions (respectively) -

$$S_{1000}(N) = \frac{1}{2} \cdot 3^N - \frac{1}{2} \quad ; \quad N > 0 \tag{6A}$$

$$P_{1000}(k) = 3^{k-1} \quad ; \quad k > 0 \ . \tag{6B}$$

On passing, note that Eq. 6A works in fact for N = 0 too [$S_{1000}(0) = 0$].

One last remark now, related to the deterministic nature of the "1000" solution, before moving on to discussing the more complicated solutions of the other far more challenging MToH puzzles. Namely, those puzzles where pre-coloring leaves the Tower with one Neutral post or two or three.

#### 3.1.7. Deterministic solution

Inspection of the "1000" MToH solution reveals that it is deterministic. That is – for a forward-moving solution, each and every move is **dictated**. In other words - there is only one way to make the next move. "Forward-moving" solution is defined as a solution where all attained Tower-States are distinct – a given Tower-State is never repeated. "Tower State" is defined as the combination of disks on the posts, including disk-size and disk-color.

For the Classical ToH, the **Optimal** (shortest Duration) solution is deterministic as well. There is only one way to make the next Optimal move.

It is interesting to note that for the Optimal solution of the Classical ToH, the recurrence relations are[7] –

$$S_{ToH}(N+1) = 2 \cdot S_{ToH}(N) + 1 \;\; ; \;\; S_{ToH}(1) = 1 \qquad (7A)$$

$$P_{ToH}(k+1) = 2 \cdot P_{ToH}(k) \;\; ; \;\; P_{ToH}(1) = 1 \qquad (7B)$$

and hence

$$S_{ToH}(N) = 2^N - 1 \;\; ; \;\; N > 0 \qquad (8A)$$

$$P_{ToH}(k) = 2^{k-1} \;\; ; \;\; k > 0 \; . \qquad (8B)$$

Thus, the Optimal solution of the Classical ToH perfectly spans base 2 (Equation 8B) and is deterministic. The solution of the **Colored** MToH perfectly spans base 3 (Equation 6B) and is deterministic.

## 3.2. The Semi-Free-C MToH and its solving Optimal Algorithms

We start by reminding ourselves of the pre-coloring configuration for the Semi-Free-C MToH-puzzle.

### 3.2.1. The pre-coloring configuration of the Semi-Free-C MToH

The pre-coloring configuration for the Semi-Free-C MToH puzzle is shown in Figure 5. Here we have one Neutral post that, during the puzzle-solving trip, may take any color.

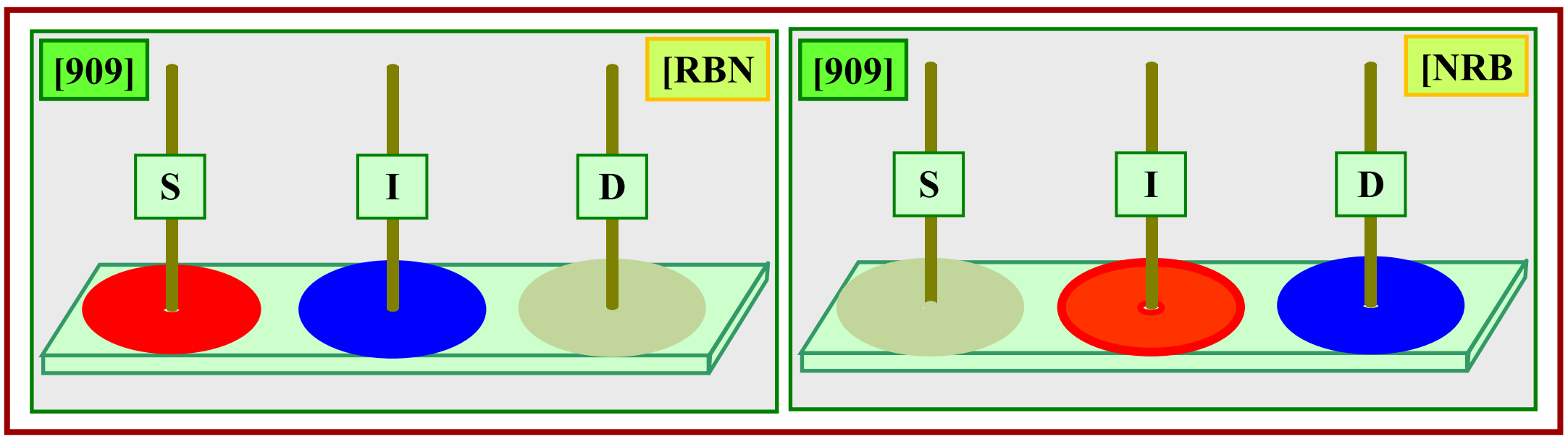

Figure 5: *The pre-coloring configurations of the Semi-Free-C MToH puzzle. The two distinct puzzle-solving Algorithms form a Time-Reversal Pair.*

Next – the Optimal "909" puzzle-solving Algorithms.

### 3.2.2. The RBN909 / NRB909 Optimal Algorithms

The Optimal Algorithms solving the Semi-Free-C MToH puzzle are listed in Table 8.

| | SRBN909 | SNRB909 |
|---|---|---|
| <br><br><br>1<br>2<br>3<br>4 | function move_SRBN909(n,s,d,i)<br>j = %N_disks + 1 - n<br>if n > 0<br>move_SRNB727(n-1,s,i,d)<br>move(j,s,d)<br>move_SRRB1000(n-1,i,s,d)<br>move_SRBB1000(n-1,s,d,i)<br>return | function move_SNRB909(n,s,d,i)<br>j = %N_disks + 1 - n<br>if n > 0<br>move_SRRB1000(n-1,s,d,i)<br>move_SRBB1000(n-1,d,i,s)<br>move(j,s,d)<br>move_SBNR727(n-1,i,d,s)<br>return |

Table 8: *The Optimal Algorithms solving the Semi-Free-N MToH-puzzle.*

The RBN pre-colored configuration is solved by the SRBN909 Algorithm listed on the left of Table 8, while the NRB pre-colored configuration is solved by the SNRB909 Algorithm listed on the right of Table 8.

The SRBN909 Algorithm and the SNRB909 Algorithm form a Time-Reversal Pair. The Tower's disk configurations

| Move # | SRBN909 | SNRB909 |
|---|---|---|
| 0 | 3 0 0 | 3 0 0 |
| 1 | 2 0 1 | 2 0 1 |
| 2 | 1 1 1 | 2 1 0 |
| 3 | 2 1 0 | 1 1 1 |
| 4 | 1 2 0 | 1 0 2 |
| 5 | 0 2 1 | 2 0 1 |
| 6 | 1 1 1 | 2 1 0 |
| 7 | 0 1 2 | 1 1 1 |
| 8 | 1 0 2 | 1 2 0 |
| 9 | 2 0 1 | 0 2 1 |
| 10 | 1 1 1 | 0 1 2 |
| 11 | 0 1 2 | 1 1 1 |
| 12 | 1 0 2 | 1 0 2 |
| 13 | 0 0 3 | 0 0 3 |

Table 9: *Number of disks on each post when solving the three-disk Semi-Free-C MToH-puzzle by each of the "909" Optimal Algorithms.*

(not including disk-size and disk color) are shown in Table 9. The two Time-Reversed columns of Table 9 are a good indication that indeed the SRBN909 Algorithm and the SNRB909 Algorithm form a Time-Reversal Pair (but they do not make a formal proof).

Note that the Tower's disk-configurations listed in Table 7 (the "1000" Pair) and those listed in Table 9 (the "909" Pair), for a three-disk puzzle, are identical. The difference between the two Algorithm-Pairs becomes evident only for N > 3.

We will now prove that each of the two Semi-Free-C MToH Solving-Algorithms listed in Table 8 is Optimal.

### 3.2.3. Proof of Optimality for the "909" Pair

The Optimality proof for both Algorithms in Table 8 is again a coupled-recursive proof. However, in the "909" case, we are looking at a more complicated case of **two Algorithm Pairs** – "909" and "727". The latter Algorithm is discussed in the next section. So to prove Optimality we show that both Pairs of Algorithms are Optimal for N = 1, we assume that **both Pairs of Algorithms** are Optimal for any n = N > 1, and we prove sequentially that each of the four Algorithms is Optimal for n = N + 1. The N + 1 part of the proof is based again on a step-by-step inspection of the MToH puzzle and using a "must" argument for each step.

We start with the SRBN909 Algorithm, on the left of Table 8.

For N = 1, both Algorithms call for one move to solve the Semi-Free-C MToH puzzle. "One move" Algorithm is obviously Optimal.

For n = N – 1 we assume Optimality of both Pairs of Algorithms (N > 2).

For n = N:

The pre-coloring configuration **for the first step** of moving N-1 disks from **S** to **I** is RNB. So the first step (line 1 in Table 8) is executed by applying the SRNB727 Algorithm to transport N-1 disks (disk 2 to disk N) from **S** to **I** using **D**. This step is Optimal by assumption. The presence of the big disk at the bottom of post **S** (was not present there in the N-1 case) makes no difference because it is big (biggest) and because it is **RED**. And we **must** perform this (Optimal) step in order to free the big disk laying on the **S**ource post and we must clear the **D**estination post.

The next step (line 2 in Table 8) is moving the big disk (disk number 1) from **S** to **D**. We **must** move any disk at least once. In this case (and in all other cases) we move the big disk exactly once which is certainly Optimal.

After performing the first two steps (lines 1 and 2 on the left of Table 8), the Tower appears Colored (for the task of moving the N-1 disks). We **must** use the "1000" Algorithms twice in the right order (lines 3 and 4 on the left of Table 8) to complete the puzzle solution. For a Colored-Tower, the "1000" Algorithms were already proved to be Optimal.

Puzzle solved. And the (conditional) proof of Optimality of the SRBN909 Algorithm ends here.

For the Time-Reversal Brother Algorithm – SNRB909 (right of Table 8), we follow the same path of proof, using similar arguments. Note that for the given pre-coloring state (NRB), we **have** to first move N-1 disks to the **BLUE** Destination post, and on to the **RED** Intermediate post and the Tower during these two moves is Colored. Only after the big disk is moved to the Destination post, the Source post becomes Neutral and we can use a "727" Algorithm, which is Optimal by assumption.

On passing, note that for the RBN (or NRB) pre-coloring combination, a less efficient yet "forward moving" solution is possible (applying the "1000" Algorithm) so that now, and certainly for progressively less restricted Towers, the solution is NOT deterministic.

The Optimality proof given above for the "909" Pair is conditional, because we still need to follow a similar Optimality proof for the "727" Pair. This is done below in section 3.3.

Before discussing the "727" Pair, we want to develop recurrence relations for the "909" Pair.

### 3.2.4. Recurrence relations for the "909" Pair

Given the "909" Solution-Algorithm Brothers (Table 8), we can extract recurrence relations for the associated number of moves.

Following the same line of arguments given in section 3.1, we find

$$S_{909}(N+1) = S_{727}(N) + 2 \cdot S_{1000}(N) + 1 \quad ; \quad S_{909}(1) = 1 \qquad (9A)$$

$$P_{909}(k+1) = P_{727}(k) + 2 \cdot P_{1000}(k) \quad ; \quad P_{909}(1) = 1 \qquad (9B)$$

The Recurrence Relations (9A and 9B) for the "909" moves involve the "727" moves. Surely the derivation of the closed-form expressions for the "909" Algorithm, must also be delayed.

We need to proceed now and analyze the "727" Algorithm Pair.

### 3.3. The Semi-Free-C MToH and its solving Optimal Algorithms

We start by reminding ourselves of the pre-coloring configuration for the Semi-Free-C MToH-puzzle.

#### 3.3.1. The pre-coloring configuration of the Semi-Free-C MToH

The pre-coloring configuration for the Semi-Free-C MToH puzzle is shown in Figure 6. Here again we have one Neutral post that, during the puzzle-solving trip, may take any color. Here we have only one pre-coloring combination but it can still be solved with two distinct Algorithms that form a Time-Reversal Pair. See next sub-section.

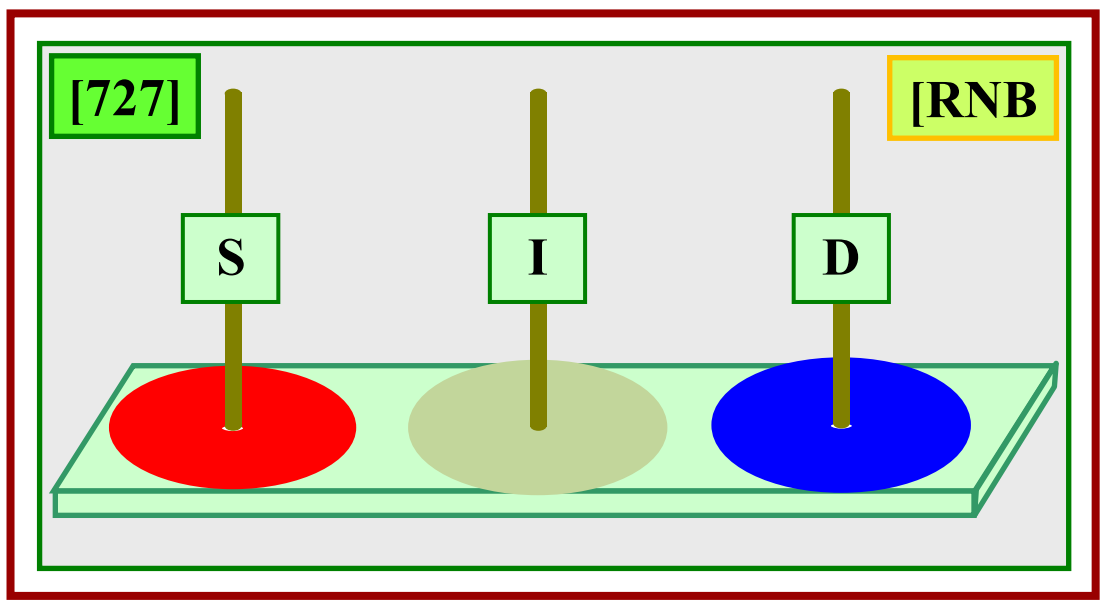

Figure 6: *The pre-coloring configuration of the Semi-Free-C MToH puzzle. Here we have only one pre-coloring combination but it can still be solved with two distinct Algorithms that form a Time-Reversal Pair.*

Next – the Optimal "727" puzzle-solving Algorithms.

### 3.3.2. The RNB727 / BNR727 Optimal Algorithms

The Optimal Algorithms solving the Semi-Free-C MToH puzzle are listed in Table 10.

| | SRNB727 | SBNR727 |
|---|---|---|
| | function move_SRNB727(n,s,d,i) | function move_SBNR727(n,s,d,i) |
| | j = %N_disks + 1 - n | j = %N_disks + 1 - n |
| | if n = 1 | if n = 1 |
| | move(j,s,d) | move(j,s,d) |
| | return | return |
| | | |
| | if n > 1 | if n > 1 |
| 1 | move_SRBN909(n-1,s,i,d) | move_SRBN909(n-2,s,i,d) |
| 2 | move(j,s,d) | move(j+1,s,d) |
| 3 | move_SRRB1000(n-2,i,s,d) | move_SNRB909(n-2,i,s,d) |
| 4 | move_SRBB1000(n-2,s,d,i) | move(j+1,d,i) |
| 5 | move(j+1,i,s) | move_SRRB1000(n-2,s,d,i) |
| 6 | move_SRBN909(n-2,d,i,s) | move_SRBB1000(n-2,d,i,s) |
| 7 | move(j+1,s,d) | move(j,s,d) |
| 8 | move_SNRB909(n-2,i,d,s) | move_SNRB909(n-1,i,d,s) |
| | return | return |

Table 10: *The Optimal Algorithms solving the Semi-Free-C MToH-puzzle.*

The RNB pre-colored configuration is solved by **both** the SRNB727 Algorithm and the SBNR727 Algorithm listed in Table 10.

The SRNB727 Algorithm and the SBNR727 Algorithm form a Time-Reversal Pair. The designation of the latter (could be somewhat confusing) signifies Time-Reversal of the former, NOT a pre-coloring combination on its own.

The Tower's disk-configurations for three disks (not including disk-size and disk color) during execution of the "727" Algorithm are shown in Table 11. The two Time-Reversed columns of Table 11 are a good indication that indeed the SRNB727 Algorithm and the SBNR727 Algorithm form a Time-Reversal Pair (but they do not make a formal proof).

Note that the "727" Algorithms solves the RNB pre-colored MToH puzzle by 11 moves only (compare with Table 7 and Table 9).

| Move # | SRNBN727 | SBNR727 |
|---|---|---|
| 0 | 3 0 0 | 3 0 0 |
| 1 | 2 0 1 | 2 1 0 |
| 2 | 1 1 1 | 1 1 1 |
| 3 | 2 1 0 | 2 0 1 |
| 4 | 1 2 0 | 2 1 0 |
| 5 | 0 2 1 | 1 1 1 |
| 6 | 1 1 1 | 1 2 0 |
| 7 | 0 1 2 | 0 2 1 |
| 8 | 1 0 2 | 0 1 2 |
| 9 | 1 1 1 | 1 1 1 |
| 10 | 0 1 2 | 1 0 2 |
| 11 | 0 0 3 | 0 0 3 |

Table 11: *Number of disks on each post when solving the three-disk Semi-Free-C MToH-puzzle by each of the "727" Optimal Algorithms.*

We will now prove that each of the two Semi-Free-C MToH Solving-Algorithms listed in Table 10 is Optimal.

### 3.3.3. Proof of Optimality for the "727" Pair

The Optimality proof for both Algorithms in Table 10 is again a coupled-recursive proof. And as already mentioned the "727" proof is coupled with the "909" proof. Here, the N = 1 case is isolated (Table 10) and the "727" Algorithms are applied to N > 1 Towers. So to prove Optimality we show that both Pairs of Algorithms are Optimal for N = 2, we assume that **both Pairs** of Algorithms are Optimal for any n = N > 2, and we prove sequentially that each Algorithm is Optimal for n = N + 1. The N + 1 part of the proof, here too, is based on step-by-step inspection of the MToH puzzle and using a "**must**" argument for each step.

We start with the SRNB727 Algorithm, on the left of Table 10.

For N = 2, both "727" Algorithms call for four moves to solve the Semi-Free-C MToH puzzle – one move for the big one and three moves for the small one. A simple set of "**must**" arguments prove that a "three and one" Algorithm is Optimal.

For n = N – 1 we assume Optimality of the "727" Pair and of the "909" Pair (N > 2).

For n = N:

The pre-coloring configuration **for the first step** of moving N-1 disks from **S** to **I,** is RBN. So the first step (line 1 in Table 10) is executed by applying the SRBN909 Algorithm to transport N-1 disks (disk 2 to disk N) from **S** to **I** using **D**. This step is Optimal by assumption. The presence of the big disk at the bottom of post **S** (was not present there in the N-1 case) makes no difference because it is big (biggest) and because it is **RED**. And we **must** perform this (Optimal) step in order to free the big disk laying on the **S**ource post and we must clear the **D**estination post.

The next step (line 2 in Table 10) is moving the big disk (disk number 1) from **S** to **D**. We **must** move any disk at least once. In this case (and in all other cases) we move the big disk exactly once which is certainly Optimal.

Steps 3 and 4 (lines 3 and 4 in Table 10) are designed to free disk number 2 (the second largest). The Optimal plan is to move disk number 2 to **S** and return the N-2 stack to the (Neutral again) **I** post, this time **RED** facing up. The Tower for these two steps is Colored so we **must** use two "1000" Algorithms (in the right order) to execute the steps.

Step number 5 (lines 5 in Table 10) is flipping the second largest disk to **RED** **S** (as planned).

Step number 6 (lines 6 in Table 10) is moving N-2 disks from **BLUE** **D** to Neutral **I** (as planned). The pre-coloring configuration is BRN which is equivalent to RBN, so we use SRBN909 for the task. The presence of the big disk at the bottom of post **D** and the second big disk at the bottom of post **S** makes no difference. This step which we **must** perform (as part of our Optimal plan) is Optimal by assumption.

Step number 7 (lines 7 in Table 10) is flipping the second largest disk to **BLUE** **D**.

The last step (lines 8 in Table 10) is moving N-2 disks from Neutral I to **BLUE** **D.** The pre-coloring configuration is NRB so we use SNRB909 for the task. We must perform this step to complete solving the puzzle and the step is Optimal by assumption.

Puzzle solved. And the (conditional) proof of Optimality of the SRNB727 Algorithm ends here.

For the Time-Reversal Brother Algorithm – SBNR727 (right of Table 10), we follow a similar path of proof, using similar arguments.

At this point the Optimality coupled recursive-proof for the "909" Pair and for the "727" Pair is completed. In subsequent Algorithms, when we run into a relevant pre-coloring configuration and execute the step using one of these four Algorithms, we **know** it is Optimal.

### 3.3.4. Recurrence relations for the "727" Pair

Given the "727" Solution-Algorithm Brothers (Table 10), we can extract recurrence relations for the associated number of moves.

Following the same line of arguments given in section 3.1, we find

$$S_{727}(N+1) = S_{909}(N) + 2\cdot S_{909}(N-1) + 2\cdot S_{1000}(N-1) + 3$$
$$S_{727}(1) = 1 \quad ; \quad S_{727}(2) = 4 \qquad (10A)$$
$$P_{727}(k+1) = P_{909}(k) + 2\cdot P_{909}(k-1) + 2\cdot P_{1000}(k-1)$$
$$P_{727}(1) = 1 \quad ; \quad P_{727}(2) = 3 \qquad (10B)$$

The Recurrence Relations (10A and 10B) for the "727" moves involve the "909" moves (see Relations 9A and 9B).

### 3.3.5. Closed-form expressions for the "909" Algorithm and for the "727" Algorithm

The recurrence relations 9A and 10A, after some algebraic manipulations, form a linear inhomogeneous recursion relations of order 3 for the $S_{727}(N)$. And similarly for $P_{727}(k)$. Once these relations are solved the relations for $S_{909}(N)$ and $P_{909}(k)$ are also determined.

The intermediate results are given by Equations 11, Equations 12, and Equations 13.

$$\lambda_1 = \sqrt[3]{1+\sqrt{26/27}} + \sqrt[3]{1-\sqrt{26/27}} \qquad (11A)$$

$$\lambda_2 = -\frac{1}{2}\cdot\lambda_1 + \frac{i}{2}\cdot\left[\sqrt[3]{\sqrt{27}+\sqrt{26}} - \sqrt[3]{\sqrt{27}-\sqrt{26}}\right] \qquad (11B)$$

$$\lambda_3 = -\frac{1}{2}\cdot\lambda_1 - \frac{i}{2}\cdot\left[\sqrt[3]{\sqrt{27}+\sqrt{26}} - \sqrt[3]{\sqrt{27}-\sqrt{26}}\right] \qquad (11C)$$

$$A_S = \frac{(7/11)\cdot\lambda_2\cdot\lambda_3 - (10/11)\cdot(\lambda_2+\lambda_3) + (19/11)}{(\lambda_2-\lambda_1)\cdot(\lambda_3-\lambda_1)} \tag{12A}$$

$$B_S = \frac{(7/11)\cdot\lambda_1\cdot\lambda_3 - (10/11)\cdot(\lambda_1+\lambda_3) + (19/11)}{(\lambda_1-\lambda_2)\cdot(\lambda_3-\lambda_2)} \tag{12B}$$

$$C_S = \frac{(7/11)\cdot\lambda_1\cdot\lambda_2 - (10/11)\cdot(\lambda_1+\lambda_2) + (19/11)}{(\lambda_2-\lambda_3)\cdot(\lambda_1-\lambda_3)} \tag{12C}$$

$$A_P = \frac{(1/11)\cdot\lambda_2\cdot\lambda_3 - (3/11)\cdot(\lambda_2+\lambda_3) + (9/11)}{(\lambda_2-\lambda_1)\cdot(\lambda_3-\lambda_1)} \tag{13A}$$

$$B_P = \frac{(1/11)\cdot\lambda_1\cdot\lambda_3 - (3/11)\cdot(\lambda_1+\lambda_3) + (9/11)}{(\lambda_1-\lambda_2)\cdot(\lambda_3-\lambda_2)} \tag{13B}$$

$$C_P = \frac{(1/11)\cdot\lambda_1\cdot\lambda_2 - (3/11)\cdot(\lambda_1+\lambda_2) + (9/11)}{(\lambda_2-\lambda_3)\cdot(\lambda_1-\lambda_3)} \tag{13C}$$

The closed-form expressions for the "909" Algorithm for N > 0 and for k > 0 are now written as –

$$S_{909}(N) = (5/11)\cdot 3^N + A_S\cdot\lambda_1^{N-1} + B_S\cdot\lambda_2^{N-1} + C_S\cdot\lambda_3^{N-1} - 1 \tag{14A}$$

$$P_{909}(k) = (10/11)\cdot 3^{k-1} + A_p\cdot\lambda_1^{k-1} + B_p\cdot\lambda_2^{k-1} + C_p\cdot\lambda_3^{k-1} \tag{14B}$$

and the closed-form expressions for the "727" Algorithm for N > 0 and for k > 0 as –

$$S_{727}(N) = (4/11)\cdot 3^N + A_S\cdot\lambda_1^{N} + B_S\cdot\lambda_2^{N} + C_S\cdot\lambda_3^{N} - 1 \tag{15A}$$

$$P_{727}(k) = (8/11)\cdot 3^{k-1} + A_p\cdot\lambda_1^{k} + B_p\cdot\lambda_2^{k} + C_p\cdot\lambda_3^{k} \tag{15B}$$

The analysis of the coupled "909" Algorithm and the "727" Algorithm is now completed.

Next - the "636" Algorithm solving the Nearly-Free MToH puzzle.

### 3.4. The Nearly-Free MToH and its solving Optimal Algorithms

We start by reminding ourselves of the pre-coloring configuration for the Nearly-Free MToH-puzzle.

#### 3.4.1. The pre-coloring configuration of the Nearly-Free MToH

The pre-coloring configuration for the Nearly-Free MToH puzzle is shown in Figure 7. Here we have two Neutral posts that, during the puzzle-solving trip, may take any color.

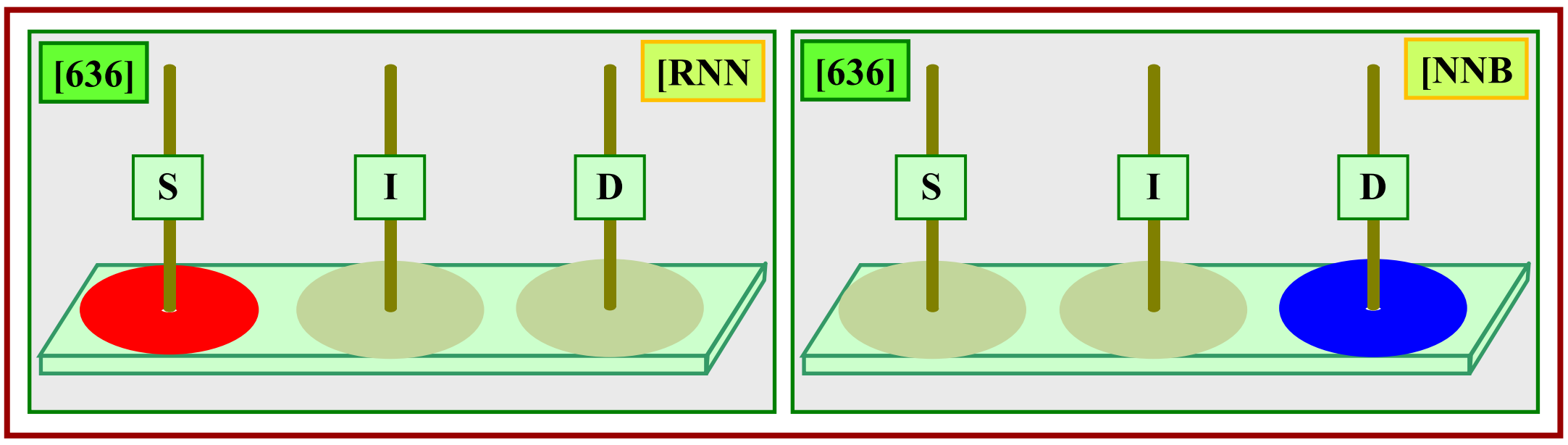

Figure 7: *The pre-coloring configuration of the Nearly-Free MToH puzzle. The Optimal Algorithms "take advantage" of the two Neutral posts and solve the puzzle with Duration of less than 64%.*

Next – the Optimal "727" puzzle-solving Algorithms.

#### 3.4.2. The RNN636 / NNB636 Optimal Algorithms

The Optimal Algorithms solving the Nearly-Free MToH puzzle are listed in Table 12. The RNN pre-colored configuration is solved by the SRNN636 Algorithm and the NNB pre-colored configuration is solved by the SNNB636 Algorithm.

**SRNN636**

```
function move_SRNN636(n,s,d,i)
  j = %N_disks + 1 - n
  if n = 1
    move(j,s,d)
    return

  if n > 1
    move_SRNN636(n-1,s,i,d)
    move(j,s,d)
    move_SRRB1000(n-2,i,s,d)
    move_SRBB1000(n-2,s,d,i)
    move(j+1,i,s)
    move_SRBN909(n-2,d,i,s)
    move(j+1,s,d)
    move_SNRB909(n-2,i,d,s)
    return
```

**SNNB636**

```
function move_SNNB636(n,s,d,i)
  j = %N_disks + 1 - n
  if n = 1
    move(j,s,d)
    return

  if n > 1
    move_SRBN909(n-2,s,i,d)
    move(j+1,s,d)
    move_SNRB909(n-2,i,s,d)
    move(j+1,d,i)
    move_SRRB1000(n-2,s,d,i)
    move_SRBB1000(n-2,d,i,s)
    move(j,s,d)
    move_SNNB636(n-1,i,d,s)
    return
```

Table 12: *The Optimal Algorithms solving the Nearly-Free MToH-puzzle.*

The Tower's disk-configurations for three disks (not including disk-size and disk color) are shown in Table 13. The two Time-Reversed columns of Table 13 are a good indication that indeed the SRNN636 Algorithm and the SNNB636 Algorithm form a Time-Reversal Pair (but they do not make a formal proof).

| Move # | SRNNN636 | SNNB636 |
|---|---|---|
| 0 | 3 0 0 | 3 0 0 |
| 1 | 2 0 1 | 2 1 0 |
| 2 | 1 1 1 | 1 1 1 |
| 3 | 2 1 0 | 2 0 1 |
| 4 | 1 2 0 | 2 1 0 |
| 5 | 0 2 1 | 1 1 1 |
| 6 | 1 1 1 | 1 2 0 |
| 7 | 0 1 2 | 0 2 1 |
| 8 | 1 0 2 | 0 1 2 |
| 9 | 1 1 1 | 1 1 1 |
| 10 | 0 1 2 | 1 0 2 |
| 11 | 0 0 3 | 0 0 3 |

Table 13: *Number of disks on each post when solving the three-disk Nearly-Free MToH-puzzle by each of the "636" Optimal Algorithms.*

We will now prove that each of the two Nearly-Free MToH Solving-Algorithms listed in Table 12 is Optimal.

### 3.4.3. Proof of Optimality for the "636" Pair

The Optimality proof for either Algorithm in Table 12 is a NON-coupled-recursive proof. Seven of the eight steps in the "636" Algorithms are done with (now proved to be) Optimal Algorithms. We just need to prove that these steps are all **necessary** for an Optimal solution (and of course sufficient to solve the puzzle).

Here again, the N = 1 case is isolated (Table 12) and the "636" Algorithms are applied to N > 1 Towers. So again, to prove Optimality of one such Algorithm, we show that it is Optimal for N = 2, we assume that the Algorithm is Optimal for any n = N > 2, and we prove that the Algorithm is Optimal for n = N + 1. The N + 1 part of the proof, here too, is based on step-by-step inspection of the MToH puzzle and using a "**must**" argument for each step.

We start with the SRNN636 Algorithm, on the left of Table 12.

For N = 2, the SRNN636 Algorithm calls for four moves to solve the Nearly-Free MToH puzzle – one move for the big one and three moves for the small one. A simple set of "**must**" arguments prove that a "three and one" Algorithm is Optimal.

For n = N – 1 we assume Optimality of the SRNN636 Algorithm (N > 3).

For n = N:

The pre-coloring configuration **for the first step** of moving N-1 disks from **S** to **I,** is RNN. So the first step (line 1 in Table 10) is executed by applying the SRNN636 Algorithm to transport N-1 disks (disk 2 to disk N) from **S** to **I** using **D**. This step is Optimal by assumption. The presence of the big disk at the bottom of post **S** (was not present there in the N-1 case) makes no difference because it is big (biggest) and because it is **RED**. And we **must** perform this (Optimal) step in order to free the big disk laying on the **S**ource post and we must clear the **D**estination post.

The next step (line 2 in Table 12) is moving the big disk (disk number 1) from **S** to **D**. We **must** move any disk at least once. In this case (and in all other cases) we move the big disk exactly once which is certainly Optimal.

From here on to puzzle solution, the line of arguments follows exactly the same line of arguments presented in the "727" proof. Because from now on the Source post is **RED** by pre-coloring and the Destination post is **BLUE** because of the presence of the big disk (**BLUE** facing up). So we are back to an RNB "727" situation.

The difference then, we now realize, between the "636" Algorithm and the "727" Algorithm is in the first step of transporting "down" N-1 disks. In the "727" it is done by "909" while in the "636" it is done by the "636" itself which is more efficient. This is why the "636" is shorter than the "727".

For the Time-Reversal Brother Algorithm – SNNB636 (right of Table 12), we follow a similar path of proof, using similar arguments.

At this point the Optimality proof for the "636" Pair is completed.

Next - recurrence relations and closed form expressions for the "636" Algorithm.

### 3.4.4. Recurrence relations for the "636" Pair

Given the "636" Solution-Algorithm Brothers (Table 12), we can extract recurrence relations for the associated number of moves.

Following the same line of arguments given in section 3.1, we find

$$S_{636}(N+1) = S_{636}(N) + 2 \cdot S_{909}(N-1) + 2 \cdot S_{1000}(N-1) + 3$$

$$S_{636}(1) = 1 \;\; ; \;\; S_{636}(2) = 4 \qquad (16A)$$

$$P_{636}(k+1) = P_{636}(k) + 2 \cdot P_{909}(k-1) + 2 \cdot P_{1000}(k-1)$$

$$P_{636}(1) = 1 \;\; ; \;\; P_{636}(2) = 3 \qquad (16B)$$

### 3.4.5. Closed-form expressions for the "636" Algorithm

Closed-form expressions for the "636" Algorithm are derived from the Recurrence-Relations 16 for N > 0 and for k > 0:

$$S_{636}(N) = (7/22) \cdot 3^N + A_S \cdot (\lambda_1 + 1) \cdot \lambda_1^{N-1} + B_S \cdot (\lambda_2 + 1) \cdot \lambda_2^{N-1}$$

$$+ C_S \cdot (\lambda_3 + 1) \cdot \lambda_3^{N-1} - (3/2) \qquad (17A)$$

$$P_{636}(k) = (7/11) \cdot 3^{k-1} + A_P \cdot (\lambda_1 + 1) \cdot \lambda_1^{k-1} + B_P \cdot (\lambda_2 + 1) \cdot \lambda_2^{k-1}$$

$$+ C_P \cdot (\lambda_3 + 1) \cdot \lambda_3^{k-1} \qquad (17B)$$

The analysis of the "636" Algorithm is now completed.

Next - the "606" Algorithm solving the Free MToH puzzle.

### 3.5. The Free MToH and its solving Optimal Algorithms

We start by reminding ourselves of the pre-coloring configuration for the Free MToH-puzzle.

#### 3.5.1. The pre-coloring configuration of the Free MToH

The pre-coloring configuration for the Free MToH puzzle is shown in Figure 8. Here we have three Neutral posts that, during the puzzle-solving trip, may take any color.

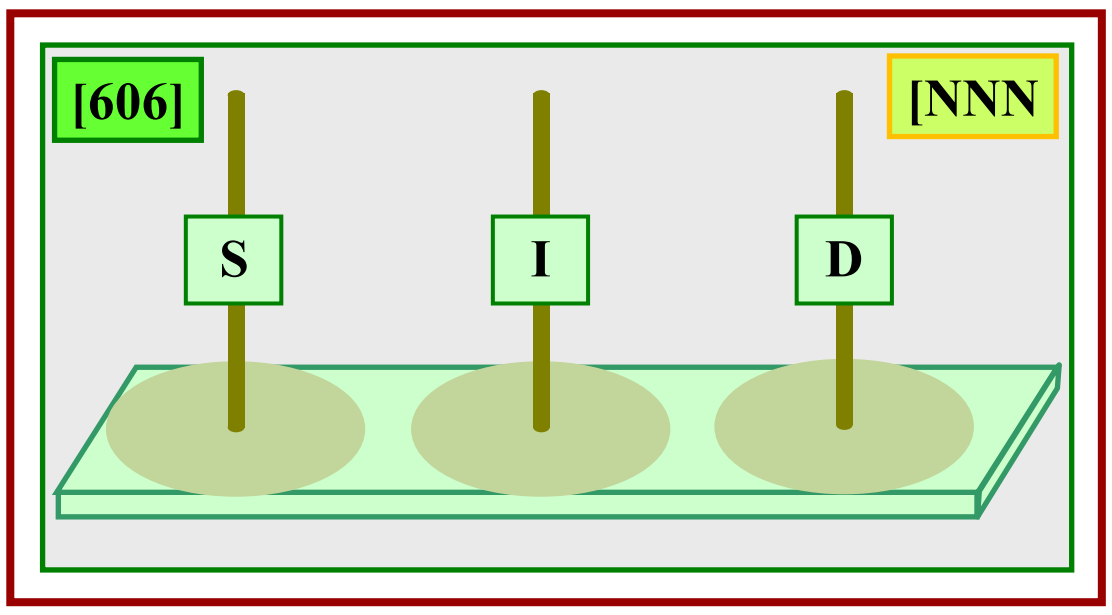

Figure 8: *The pre-coloring configuration of the Free MToH puzzle. The Optimal Algorithms "take advantage" of the three Neutral posts and solve the puzzle with Duration of less than 61%.*

Next – the Optimal "606" puzzle-solving Algorithms.

#### 3.5.2. The "606" Optimal Algorithms

The Optimal Algorithms solving the Free MToH puzzle are listed in Table 14. The "solve_up_MToH_puzzle_SNNN606" Algorithm and the "solve_down_MToH_puzzle_SNNN606" Algorithm solve the NNN pre-colored configuration. The two Algorithms form a Time-Reversal Pair.

| | solve_up_MToH_puzzle_SNNN606 | solve_down_MToH_puzzle_SNNN606 |
|---|---|---|
| | **Function** | **Function** |
| | solve_up_MToH_puzzle_SNNN606(n,s,d,i) | solve_down_MToH_puzzle_SNNN606(n,s,d,i) |
| | j = %N_disks + 1 - n | j = %N_disks + 1 - n |
| | if n > 0 | if n > 0 |
| 1 | move_SRNN636(n-1,s,i,d) | move_all_but_n_down_606(n-1,s,i,d) |
| 2 | move(j,s,d) | move(j,s,d) |
| 3 | move_all_but_n_up_606(n-1,i,d,s) | move_SNNB636(n-1,i,d,s) |
| | **function** | **function** |
| | move_all_but_n_up_606(n,s,d,i) | move_all_but_n_down_606(n,s,d,i) |
| | j = %N_disks + 1 - n | j = %N_disks + 1 - n |
| | if n > 0 | if n > 0 |
| 4 | move_SRRB1000(n-1,s,i,d) | move_SRNN636(n-1,s,d,i) |
| 5 | move_SRBB1000(n-1,i,d,s) | move(j,s,i) |
| 6 | move(j,s,i) | move_SNRB909(n-1,d,s,i) |
| 7 | move_SRBN909(n-1,d,s,i) | move(j,i,d) |
| 8 | move(j,i,d) | move_SRRB1000(n-1,s,i,d) |
| 9 | move_SNNB636(n-1,s,d,i) | move_SRBB1000(n-1,i,d,s) |

Table 14: *The Optimal Algorithms solving the Free MToH-puzzle.*

The Tower's disk-configurations for three disks (not including disk-size and disk color) are shown in Table 15. The two Time-Reversed columns of Table 11 are a good indication that indeed the "up-606" Algorithm and the "down-606" Algorithm form a Time-Reversal Pair (but they do not make a formal proof).

| Move # | solve_up_MToH_puzzle_SNNN606 | solve_down_MToH_puzzle_SNNN606 |
|---|---|---|
| 0 | 3 0 0 | 3 0 0 |
| 1 | 2 0 1 | 2 1 0 |
| 2 | 1 1 1 | 1 1 1 |
| 3 | 2 1 0 | 2 0 1 |
| 4 | 1 2 0 | 2 1 0 |
| 5 | 0 2 1 | 1 1 1 |
| 6 | 1 1 1 | 1 2 0 |
| 7 | 0 1 2 | 0 2 1 |
| 8 | 1 0 2 | 0 1 2 |
| 9 | 1 1 1 | 1 1 1 |
| 10 | 0 1 2 | 1 0 2 |
| 11 | 0 0 3 | 0 0 3 |

Table 15: *Number of disks on each post when solving the three-disk Free MToH-puzzle by each of the "606" Optimal Algorithms.*

We will now prove that each of the two Free MToH Solving-Algorithms listed in Table 14 is Optimal.

### 3.5.3. Proof of Optimality for the "606" Pair

The Optimality proof for both Algorithms in Table 14 is a NOT a recursive proof. The two steps in the first function and the six steps in the second function of each Algorithm are **all** done with (now proved to be) Optimal Algorithms. And the two Algorithms are independent of each other. We just need to prove that the steps for each Algorithm are all **necessary** for an Optimal solution (and of course sufficient to solve the puzzle).

The "up-606" Algorithm (left of Table 14) is made of two functions – "solve_up_MToH_puzzle_SNNN606" that calls for the services of "move_all_but_n_up_606".

For N = 1, only the "move(1,s,d)" step in the first function (line number 2 in Table 14) is executed. The puzzle is solved with one move which is Optimal.

For N > 1:

The first step of the first function calls for moving N-1 disk from **S** to **I**. The pre-coloring configuration **for the step** is RNN. We execute this step with a "636" Algorithm (line 1 in Table 14) which, we now **know,** is Optimal. We **must** perform this step in order to free the big disk laying on the **S**ource post and we must clear the **D**estination post.

The next step (line 2 in Table 14) is moving the big disk (disk number 1) from **S** to **D**. We **must** move the big disk to the Destination post. And one-count move is certainly Optimal.

And now the services of the second function are called for (to move N-1 disks from the **BLUE** Intermediate post to the **BLUE** Destination post).

The first two moves of the second function (line 4 and line 5 in Table 14) clear the way for the second big disk on post **I**. We **must** perform these moves and since the pre-coloring of the Tower during these two steps is NBB which is equivalent to RBB, we use the suitable (Optimal in this case) 1000 Algorithms.

Next, disk number 2 is moved from **I** to **S** (line 6 in Table 14). Note that "move_all_but_n_up_606(n,s,d,i)" is called with n = N-1, its Source post is the original **I** and its Destination post is the original **S**.

Now, the Destination post is cleared by moving N-2 disks (N-1 for "move_all_but_n_up_606") from original **D** to original **I** (line 7 in Table 14). The pre-coloring configuration **for this step** is (equivalent to) RBN so we use the Optimal "909" Algorithm to execute this **necessary** step.

The second largest disk is moved from the (original) Source post to the (original) Destination post (line 8 in Table 14).

Finally, we **must** move the N-2 disks (N-1 for "move_all_but_n_up_606") from original **I** to original **D** (line 9 in Table 14). The pre-coloring configuration **for this step** is BBN so we use the Optimal "636" Algorithm to execute this step.

Puzzle solved.

Optimality proof for the "solve_up_MToH_puzzle_SNNN606" Algorithm (left of Table 14) ends here.

The "solve_down_MToH_puzzle_SNNN606" (right of Table 14) is an independent Algorithm, actually unnecessary and is listed here only to show that a second, yet equivalent, solution-route exists. Multiplicity of solution-routes is discussed further in section 4 below.

Optimality proof for the "606-down" solution follows a course very similar to the course for the Optimality proof for the "606-up" solution.

At this point the independent Optimality proof for each of the "606" Pair is completed.

And at this point we know that a **Duration-Limit of 606/1000** (exactly 20/33 – see below) **is the minimum Duration-Limit for solving the Free Magnetic Tower of Hanoi**.

A "movie" showing how a five disk MToH-puzzle is solved by a programmed "606" Algorithm in 83 moves is now on the web[8].

Next - recurrence relations and closed form expressions for the "606" Algorithm.

### 3.5.4. Recurrence relations for the "606" Algorithm

Given the "606" Solution-Algorithm (either side of Table 14), we can extract recurrence relations for the associated number of moves:

$$S_{606}(N+1) = S_{636}(N) + S_{636}(N-1) + S_{909}(N-1) + 2 \cdot S_{1000}(N-1) + 3$$
$$S_{606}(1) = 1 \quad ; \quad S_{606}(2) = 4$$

(18A)

$$P_{606}(k+1) = P_{636}(k) + P_{636}(k-1) + P_{909}(k-1) + 2 \cdot P_{1000}(k-1)$$
$$P_{606}(1) = 1 \quad ; \quad P_{606}(2) = 3 \qquad (18B)$$

### 3.5.5. Closed-form expressions for the "606" Pair

Closed-form expressions for the "606" Algorithm are derived from the Recurrence-Relations 18, using already-determined closed-form expressions. The end results are given be Equations 19:

$$
\begin{aligned}
S_{606}(N) = \left(10/33\right)\cdot 3^{N} + \left(1/2\right)\cdot A_S \cdot \left(\lambda_1+1\right)^2 \cdot \lambda_1^{N-1} \\
+\left(1/2\right)\cdot B_S \cdot \left(\lambda_2+1\right)^2 \cdot \lambda_2^{N-1} \\
+\left(1/2\right)\cdot C_S \cdot \left(\lambda_3+1\right)^2 \cdot \lambda_3^{N-1} - 2 \quad ; \quad N>0
\end{aligned}
\tag{19A}
$$

$$
P_{606}(1) = 1 \tag{19B1}
$$

$$
\begin{aligned}
P_{606}(k) = \left(20/33\right)\cdot 3^{k-1} + \left(1/2\right)\cdot A_P \cdot \left(\lambda_1+1\right)^2 \cdot \lambda_1^{k-1} \\
+\left(1/2\right)\cdot B_P \cdot \left(\lambda_2+1\right)^2 \cdot \lambda_2^{k-1} \\
+\left(1/2\right)\cdot C_P \cdot \left(\lambda_3+1\right)^2 \cdot \lambda_3^{k-1} \quad ; \quad k>1
\end{aligned}
\tag{19B2}
$$

The analysis of all five Optimal Solution-Algorithms ends here. We now know how to Optimally solve the MToH puzzle with **any** pre-coloring configuration. Particularly the original and "natural" Free (NNN) configuration.

Before concluding, we want to briefly discuss multiplicity of solution-routes.

## 4. Multiplicity of solution-routes

When solving a two-disk MToH puzzle, we quickly realize that two equivalent solution options exist – move the small disk once, move the big disk once and move the small disk twice. Or - move the small disk twice, move the big disk once and move the small disk once more. These two solution "strategies" show up for a larger number of disks too. But clearly, **in terms of number of moves, they are equivalent**.

For a three-disk Free MToH puzzle, four different solution routes exist (not explicitly presented in this paper). All four consist of (1,3,7) moves, all are equivalent, and all are Optimal. They differ only in the order in which the disks are transported.

| | 1 | 2 | 3 | 4 | 5 | 6 | 7 | 8 |
|---|---|---|---|---|---|---|---|---|
| **Step #** | **U11** | **U21** | **U12** | **U22** | **D11** | **D21** | **D12** | **D22** |
| **0** | **7 0 0** | **7 0 0** | **7 0 0** | **7 0 0** | **7 0 0** | **7 0 0** | **7 0 0** | **7 0 0** |
| **xxx** | **xxx** | **xxx** | **xxx** | **xxx** | **xxx** | **xxx** | **xxx** | **xxx** |
| **58** | **304** | **304** | **304** | **304** | **340** | **340** | **340** | **340** |
| **59** | **313** | **403** | **313** | **403** | **331** | **430** | **331** | **430** |
| **60** | **412** | **313** | **412** | **313** | **421** | **331** | **421** | **331** |
| **61** | **403** | **412** | **403** | **412** | **430** | **421** | **430** | **421** |
| **62** | **502** | **502** | **502** | **502** | **520** | **520** | **520** | **520** |
| **xxx** | **xxx** | **xxx** | **xxx** | **xxx** | **xxx** | **xxx** | **xxx** | **xxx** |
| **81** | **403** | **403** | **403** | **403** | **430** | **430** | **430** | **430** |
| **82** | **313** | **313** | **304** | **304** | **331** | **331** | **340** | **340** |
| **83** | **214** | **214** | **313** | **313** | **241** | **241** | **331** | **331** |
| **84** | **304** | **304** | **214** | **214** | **340** | **340** | **241** | **241** |
| **85** | **205** | **205** | **205** | **205** | **250** | **250** | **250** | **250** |
| **xxx** | **xxx** | **xxx** | **xxx** | **xxx** | **xxx** | **xxx** | **xxx** | **xxx** |
| **687** | **0 0 7** | **0 0 7** | **0 0 7** | **0 0 7** | **0 0 7** | **0 0 7** | **0 0 7** | **0 0 7** |

Table 16: *Eight equivalent Optimal-routes solving a seven-disk Free MToH-puzzle.*

With increased number of disks, the number of solution routes increases as well. Specifically, each time the "727" Algorithm is called, two options exist which, at the end of their execution, yield exactly the same result. These are SRNB727 and its Time Reversal Brother SBNR727 (see Table 10). Either one of the options works at the "727" pre-coloring configuration because an RNB pre-coloring configuration is in fact **identical** to a BNR pre-coloring configuration (see Figure 1). An example of eight different routes, all (Optimally) solving a seven-disk Free MToH puzzle, is given by Table 16.

## 5. Concluding remarks

The Magnetic Tower of Hanoi puzzle (or "set of puzzles") is a "colorful" extension of the Classical Tower of Hanoi. It appeals to the public, as was evident following a "Gathering for Gardner nine (2010)" talk[9], and as is manifested now in two MToH-applets voluntarily created and uploaded onto the web[5,6]. What's more - the MToH puzzle, as demonstrated throughout these pages, forms a basis for what I view as a most elegant mathematical analysis.

In this paper we have "distilled" five independent sister-puzzles and showed how to Optimally solve each of them. For each sister-puzzle we presented move-recurrence-relations, derived closed form expressions for the total (minimum) number of disk-moves required to solve the puzzle, and derived closed form expressions for the number of moves each disk is making during execution of such solution.

Analyzing the game, we find multiple (equivalent) solution-routes. We also have found several non-optimal solution Algorithms (see section 2.1 above and see reference [2]). Overall then, this wondering-allowing game presents challenges to every player.

Undoubtedly, other game extensions are "around the corner" – more-than-three posts, different move-rules, different end-state definitions, etc. Would one such extension produce a more appealing puzzle? Would one future variation lend itself to a more "colorful" analysis?

**Acknowledgements**

- Writing of this paper was triggered, and in a way guided, by the professional and constructive criticism written by a Referee of a Mathematical Journal, who reviewed a previous (and very different) version. I wish to thank the Referee for his thoughtful comments.
- Many thanks to **Professor Ruth Lawrence-Naimark** - Associate Professor of mathematics at the Einstein Institute of Mathematics, Hebrew University of Jerusalem, for **deriving the closed-form expressions** for the number of disk-moves from the respective recurrence relations. I am truly grateful to Professor Lawrence for the mathematics and no-less for her cooperative attitude.